\documentclass[leqno]{conm-p-l}
\usepackage{bbm}
\usepackage{epsf}

\newcommand{\Aff}{{\mathsf{Aff}}}
\newcommand{\A}{{\mathbbm A}}
\newcommand{\Aut}{{\mathsf{Aut}}}
\newcommand{\DD}{{\mathsf{Def}}}
\newcommand{\Def}{{\mathfrak{D}}}
\newcommand{\bigdc}{{\DD_{\mathrm c}(S,x_0)}}
\newcommand{\bigdcp}{{\DD'_{\mathrm c}(S,x_0)}}
\newcommand{\bigd}{{\DD(S,x_0)}}
\newcommand{\bigdp}{{\DD'(S,x_0)}}

\newcommand{\ddx}{\frac{\partial}{\partial x}}
\newcommand{\ddy}{\frac{\partial}{\partial y}}

\newcommand{\dev}{{\mathsf{dev}}}
\newcommand{\Diffb}{{\dif(S,x_0)}}
\newcommand{\Diff}{{\dif(S)}}
\newcommand{\Diffob}{{\dif_0(S,x_0)}}
\newcommand{\DiffOb}{{\dif^0(S,x_0)}}
\newcommand{\Diffo}{{\dif_0(S)}}

\newcommand{\dif}{{\mathsf{Diff}}}

\newcommand{\DT}{{\DD_{\mathrm c}(T^2)}}

\newcommand{\GL}{{\mathsf{GL}}}
\newcommand{\GLtr}{{\mathsf{GL}(2,\R)}}
\newcommand{\GLtrp}{{\mathsf{GL}_+(2,\R)}}

\newcommand{\h}{{\mathsf{H}^2}}

\newcommand{\hol}{{\mathsf{hol}}}
\newcommand{\Hol}{{\mathsf{Hol}}}
\newcommand{\Hom}{{\mathsf{Hom}}}
\newcommand{\hpg}{{\Hom(\pi,\Aff(n))}}
\newcommand{\hpgp}{{\Hom_{\mathrm{p,S}}(\pi,\Aff(n))}}

\newcommand{\Nrm}{{\mathsf{Norm}_+(G_1)}}

\newcommand{\pp}{{\mathsf{p}}}

\newcommand{\R}{{\mathbb R}}

\newcommand{\SLtz}{{\mathsf{SL}(2,\Z)}}
\newcommand{\tM}{{\tilde{M}}}
\newcommand{\tS}{{\tilde{S}}}

\newcommand{\Vect}{{\mathsf{Vect}(M)}}
\newcommand{\vk}{{\mathsf{k}}}
\newcommand{\vp}{{\mathsf{p}}}
\newcommand{\vq}{{\mathsf{q}}}
\newcommand{\vS}{{\mathsf{S}}}

\newcommand{\xmapsto}[1]{|\hspace{-5pt}\xrightarrow{~{#1}~}}

\newcommand{\Z}{{\mathbb Z}}

\newtheorem{thm}{Theorem}[section]
\newtheorem{proposition}[thm]{Proposition}
\newtheorem{lemma}[thm]{Lemma}

\numberwithin{equation}{section}

\begin{document}

\title[Is The Deformation Space Smooth?]{Is the deformation space of complete affine structures on the 2-torus smooth?}

\author{Oliver Baues}
\address{ 
    Departement Mathematik, 
    ETH-Zentrum,  
    R\"amistrasse 101 \\
    CH-8092 Z\" urich, SWITZERLAND }
\email{ oliver@math.ethz.ch}

\author{William M. Goldman}
\address{ Mathematics Department \\
    University of Maryland, College Park, MD  20742 USA  }
\email{ wmg@math.umd.edu }

\thanks{Goldman gratefully
acknowledges partial support from National Science Foundation grant
DMS-0103889.  This research was performed at the S\'eminaire Sud-Rhodanien at
the Centre international de recherches en math\'ematiques in Luminy in November
2001, and the workshop ``Discrete Groups and Geometric Structures, with
Applications'' held at the Katolieke Universiteit Leuven Campus Kortrijk, in
Kortrijk, Belgium in March 2002.  We are grateful for the hospitality afforded
by these two meetings.}

\dedicatory{Dedicated to Alberto Verjovsky on his sixtieth birthday}


\begin{abstract} 
Periods of parallel exterior forms define natural coordinates on the
deformation space of complete affine structures on the two-torus. 
These coordinates define a differentiable 
structure on this deformation space, under which it is diffeomorphic
to $\R^2$. 
The action of the mapping class group of $T^2$ is equivalent in these
coordinates with the standard linear action of $\SLtz$ on $\R^2$.
\end{abstract}

\maketitle

\section{Introduction}
Conformal structures on the 2-torus $T^2$ (elliptic curves) are
classified by a moduli space, which is the quotient of the
upper half-plane $\h$ by the action of the modular group $\SLtz$ by
linear fractional transformations.
In other words, equivalence classes of elliptic curves correspond
to orbits of $\SLtz$ on $\h$.
Since $\h/\SLtz$ is not a smooth manifold, 
it is often easier to study properties of elliptic 
curves in terms of the {\em action\/} of $\SLtz$ on the 
smooth manifold $\h$.

This note concerns the analogous question for {\em complete affine
structures\/} on $T^2$.  A {\em complete affine structure\/} on a
manifold $M$ is a representation of $M$ as the quotient of affine
space $\A^n$ by a discrete group of affine transformations.
Equivalence classes of such structures on $T^2$ identify with 
the orbits of $\SLtz$ on $\R^2$ for the standard linear action of $\SLtz$ on
$\R^2$. However, unlike the action on $\h$, this action is {\em not
proper\/} and the quotient space is even more badly behaved.  In
particular it fails to be Hausdorff and since the action of 
$\SLtz$ on $\R^2$ is ergodic \cite[2.2.9]{Zimmer}, it also 
enjoys a nonstandard Borel structure \cite[2.1.16]{Zimmer},
\cite[\S 3]{Arveson}. 

The coordinates in $\h$, parametrizing elliptic curves,
are defined by the periods
of holomorphic 1-forms. Analogously, we use periods of parallel
1-forms to define coordinates on the deformation space of 
complete affine structures on $T^2$.  
These coordinates define a natural differentiable 
structure on this deformation space, under which it is diffeomorphic
to $\R^2$. (Using a different approach 
Baues~\cite{BauesG} noticed that this deformation space
is homeomorphic to $\R^2$.)

Whereas every smooth deformation of elliptic curves 
is induced by a \emph{smooth} map into $\h$,  
the analogous result fails for the deformation space of 
complete affine structures on $T^2$ with respect to the
coordinates given by the periods of parallel 1-forms.
We provide an example of a smooth one-parameter deformation of an 
Euclidean 2-torus, in fact an affine fibration, which 
corresponds to  a continuous curve in the deformation space
which is \emph{not} smooth in the origin.

The complete affine structures on $T^2$ were classified by
Kuiper~\cite{Kuiper}. In this paper we describe the {\em moduli
space\/} of such structures. As is often the case in moduli problems,
the moduli space can be made more tractable by introducing extra
structure --- in this case a {\em marking}, and the moduli problem
for marked structures is simpler. Equivalence classes of marked
structures form a space homeomorphic to the plane $\R^2$, which we
call the {\em deformation space.\/} The moduli space --- consisting of
equivalence classes of structures (without marking) --- is then the
orbit space of $\R^2$ by the standard linear action of $\SLtz$.

That 
the deformation space of complete affine structures on
$T^2$ is even Hausdorff is in itself surprising. (The larger
deformation space of (not necessarily complete) affine structures 
on $T^2$ is {\em not\/} Hausdorff.) 

Moreover, in our setup 
the group actions defining 
the moduli problem are not reductive, and hence the usual 
techniques of geometric invariant theory do not apply. 
Indeed, the affine structures discussed here all come 
from affine structures on abelian Lie groups which are invariant 
under multiplication, and the resulting automorphisms of the 
affine structure form a unipotent group.  \\

Benzecri \cite{Ben} (see also Milnor~\cite{Milnor})
proved that the only compact orientable surface which admits 
affine structures is the two-torus $T^2$.  The classification
of affine structures on $T^2$ up to affine diffeomorphism, 
initiated by Kuiper \cite{Kuiper}, was completed independently by 
Furness-Arrowsmith \cite{FurArr} and Nagano-Yagi \cite{NaganoYagi}.

According to Kuiper%~\cite{Kuiper}
, the complete affine structures on the 2-torus fall into two types.
Quotients of $\A^2$ by translations we call {\em Euclidean structures,\/}
since they admit compatible Euclidean (or flat Riemannian) metrics. 
The Levi-Civita connection of a compatible metric 
coincides with the flat connection associated to the affine
structure. The Euclidean structures fall into a single affine 
equivalence class, which corresponds to the origin in $\R^2$. 
The other structures are more
exotic, and do not correspond to Levi-Civita connections of a 
(possibly indefinite) Riemannian metric. 

The non-Riemannian structures can be characterized as follows.
Let $M$ be a marked non-Riemannian affine 2-torus. 
There is a parallel area form $\omega_M$ which we normalize to
have area $1$. Then there exists a vector field $\xi_M$, 
polynomial of degree $1$, such that the covariant derivative 
$\zeta_M = \nabla_{\xi_M}\xi_M$ 
is a nonzero parallel vector field. Since parallel vector fields
on $M$ are determined up to multiplication by a nonzero scalar,
we can normalize $\zeta_M$ by requiring that
\begin{equation*}
\omega_M(\xi_M,\zeta_M) = 1. 
\end{equation*}
Then the interior product
\begin{equation*}
\eta_M\, := \, \iota_{\zeta_M}\,(\omega_M)  
\end{equation*}
is a parallel $1$-form on $M$, uniquely determined by
the affine structure of $M$. The cohomology class 
\begin{equation*}
[\eta_M] \in H^1(T^2;\R)
\end{equation*}
determines the equivalence class of the marked affine manifold $M$.
As non-Euclidean structures converge to the Euclidean structure,
this cohomology class approaches zero, and the period mapping extends
to a homeomorphism of the deformation space of {\em all\/} complete
marked affine structures on $T^2$ with $\R^2$.

The case when the periods are {\em rationally related,\/} that
is when $[\eta_M]$ lies in $\R\cdot H^1(T^2;\Z)$, is particularly notable.
In this case the parallel vector field $\zeta_M$ has closed trajectories,
and the homology class of a trajectory is Poincar\'e dual to $[\eta_M]$.
Equivalently, $\eta_M$ integrates to an affine fibration 
$M\longrightarrow S^1$, where $S^1$ is given a complete affine structure.\\

For every complete affine structure on $T^2$ there is a commutative group
operation on $T^2$ for which multiplication is affine.  Thus every
complete affine 2-torus $M$ is a group object in the category of affine
manifolds. Such objects correspond to associative algebras, providing
an alternative approach to the classification. From this viewpoint, a
marked complete affine structure on $T^2$ corresponds to a
commutative associative nilpotent $\R$-algebra, with a chosen basis
corresponding to the marking. 

Similarly Matsushima \cite{Matsushima}
showed that the set of all complex affine structures on a 
fixed {\em complex torus\/} forms an affine algebraic variety,  
by identifying each structure with the associative 
algebra which is formed by the holomorphic vector fields of the 
torus under the covariant differentiation of the affine structure. 

As a related fact the deformation spaces of complete affine structures
on (also higher dimensional) real tori are homeomorphic to 
a real algebraic set \cite{BauesV}. \\

In figure 1-9  we show  various tilings of the plane 
which are obtained from affine group actions
corresponding to different values of the period class. 
The software used for creating these pictures is available
from W.\ Goldman at the Experimental Geometry Lab of the
University of Maryland.

%The Mathematica notebook used for creating these pictures 
%is available from 
%\begin{quote}
%\texttt{http://www.math.umd.edu/$\sim$res/EGL/AffineT2.nb}.
%\end{quote}
%An interactive Java applet for investigating these structures 
%can be found at 
%\begin{quote}
%\texttt{http://www.math.umd.edu/$\sim$lidador/Affine/}.
%\end{quote}
%%%%%%%%%%%%%%%%%%%%%%%%%%%%%%%%%%%%%%%%%%%%%%%%%%%%%%%%%%%%%%%%%%%%%%%%
% \input pictures
% pictures.tex
% picture section

\clearpage
\bigskip
\begin{figure}[htp]
\centerline{\epsfxsize=3in \epsfbox{ 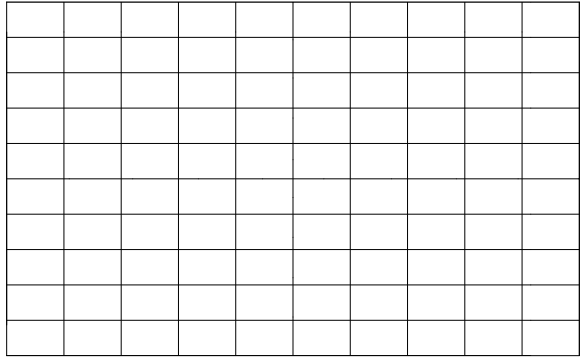}}
\caption{Period Class = (0,0)}
\label{fig:000}
\end{figure}

\bigskip
\begin{figure}[htp]
\centerline{\epsfxsize=3in \epsfbox{ 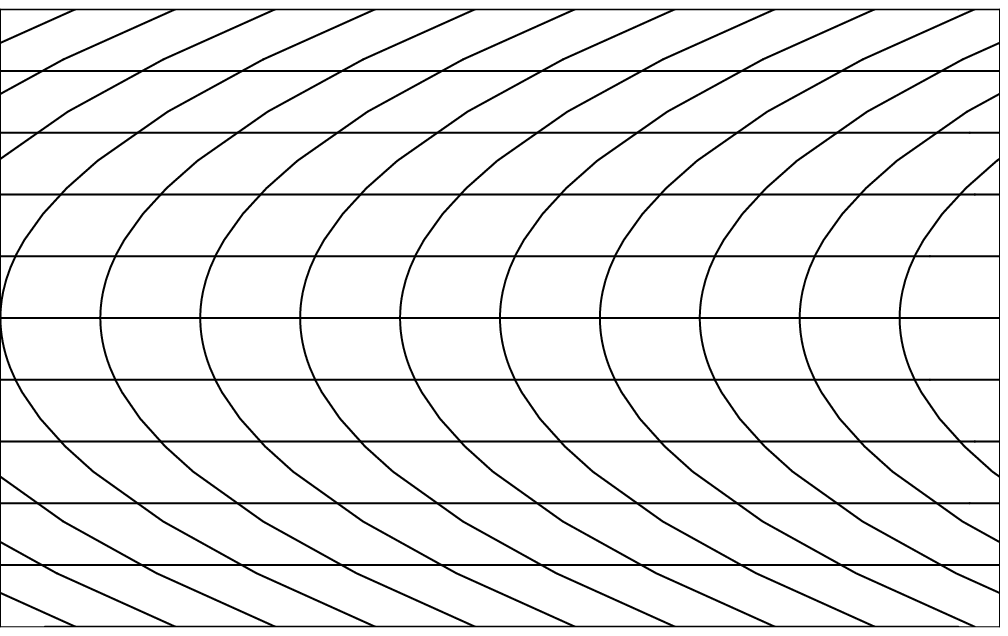}}
\caption{Period Class = (0.7,0)}
\label{fig:p00}
\end{figure}

\bigskip
\begin{figure}[htp]
\centerline{\epsfxsize=3in \epsfbox{ 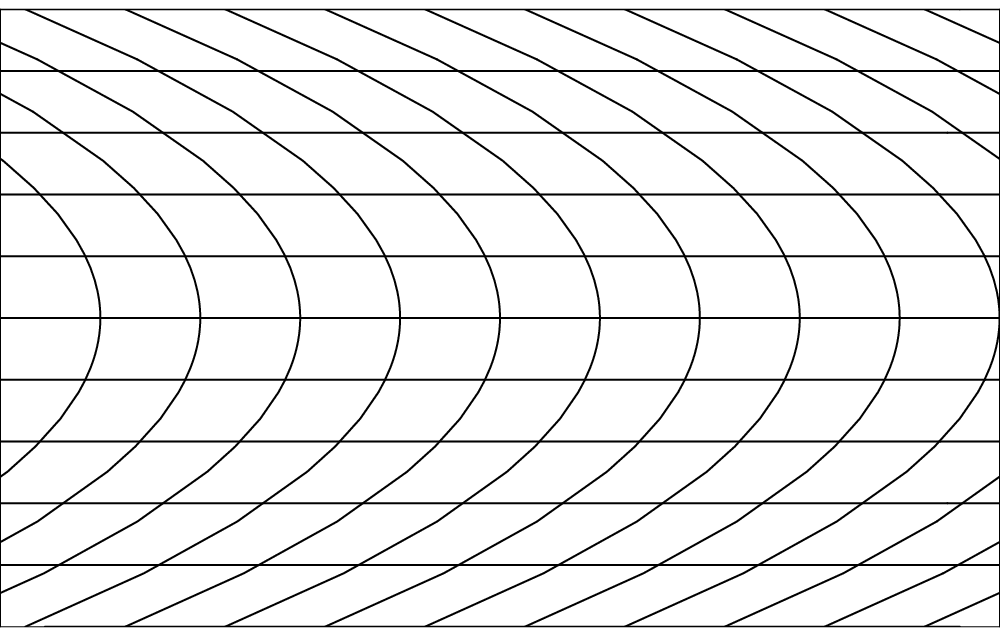}}
\caption{Period Class = (-0.7,0)}
\label{fig:m00}
\end{figure}
\newpage

\bigskip
\begin{figure}[htp]
\centerline{\epsfxsize=3in \epsfbox{ 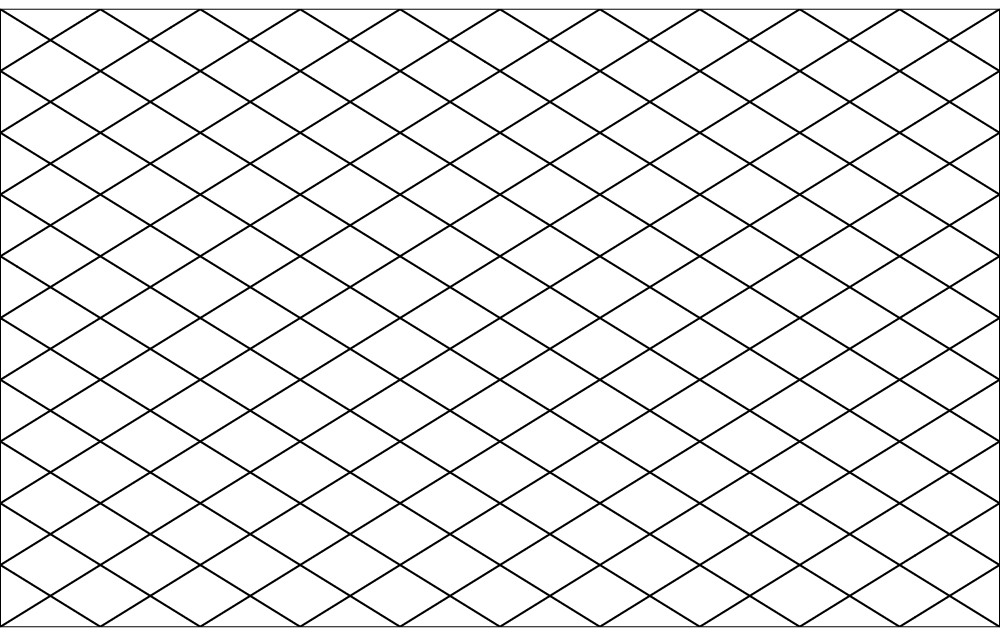}}
\caption{Period Class = (0,0)}
\label{fig:0m1}
\end{figure}

\bigskip
\begin{figure}[htp]
\centerline{\epsfxsize=3in \epsfbox{ 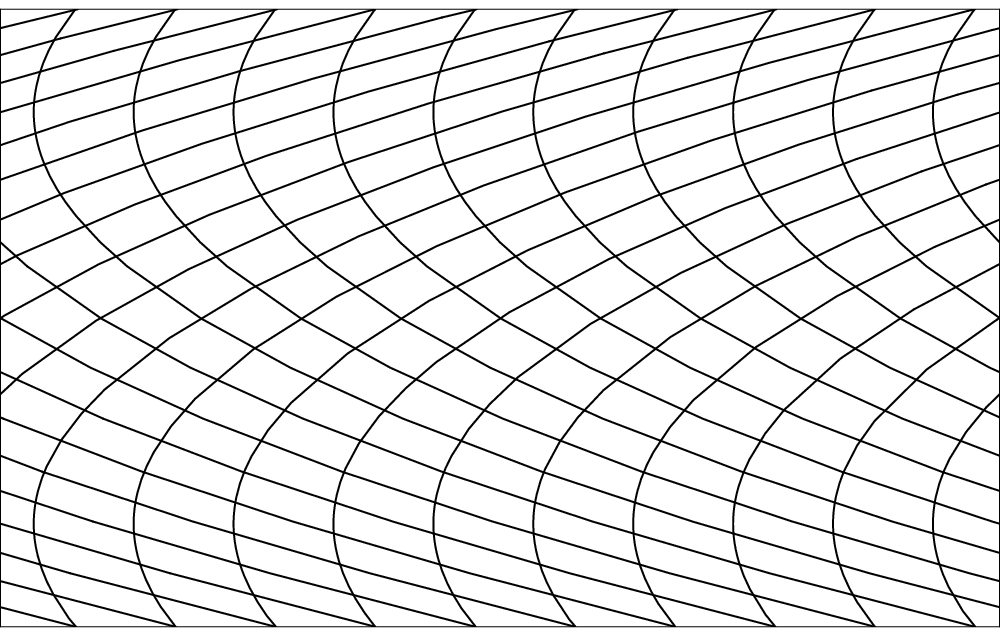}}
\caption{Period Class = (0.5,0.5)}
\label{fig:pm1}
\end{figure}

\bigskip
\begin{figure}[htp]
\centerline{\epsfxsize=3in \epsfbox{ 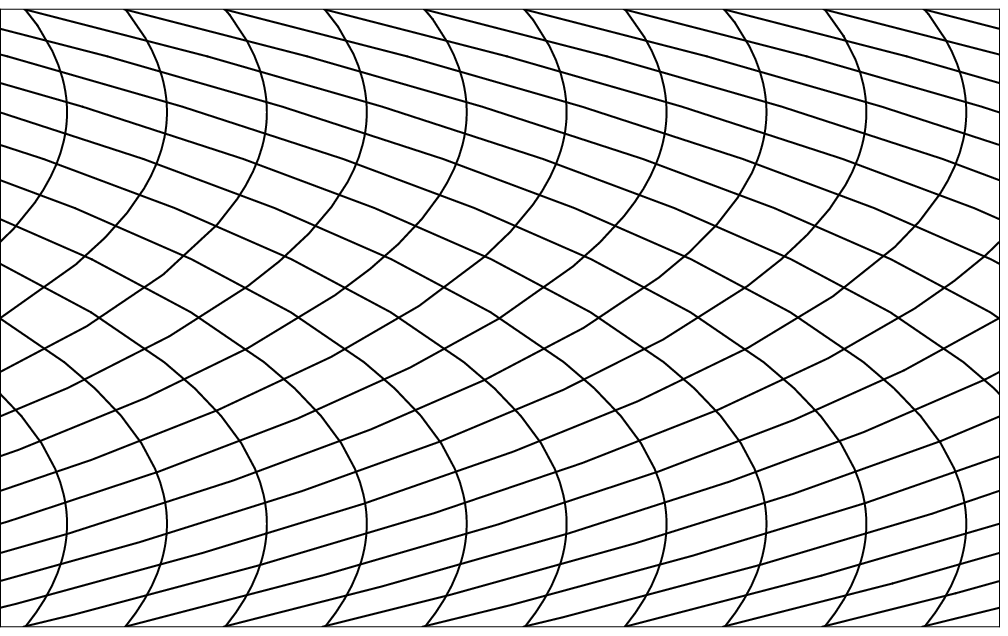}}
\caption{Period Class = (-0.5,-0.5)}
\label{fig:mm1}
\end{figure}

%%%%%%%%%%%%%%%%%%%%%%%%%%%%%%
\bigskip
\begin{figure}[htp]
\centerline{\epsfxsize=3in \epsfbox{ 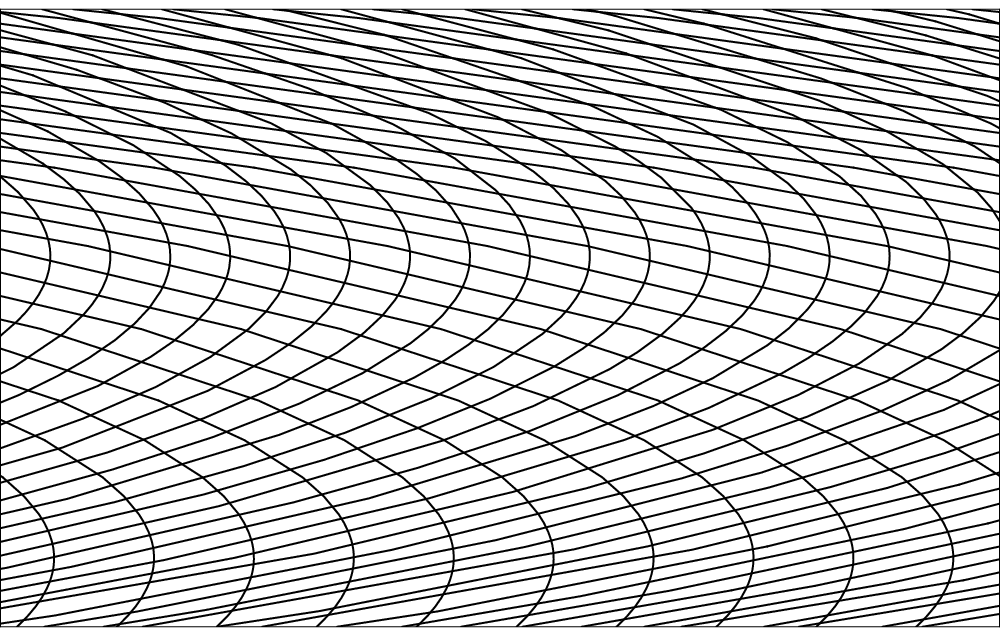}}
\caption{Period Class = (0.7,1.2)}
\label{fig:np1}
\end{figure}

\bigskip
\begin{figure}[htp]
\centerline{\epsfxsize=3in \epsfbox{ 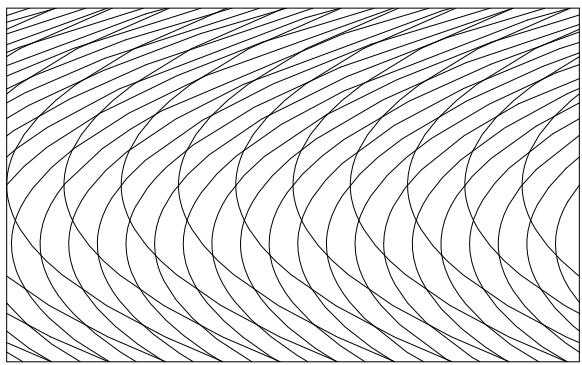}}
\caption{Period Class = (-1.7,-1.9)}
\label{fig:np2}
\end{figure}

\bigskip
\begin{figure}[htp]
\centerline{\epsfxsize=3in \epsfbox{ 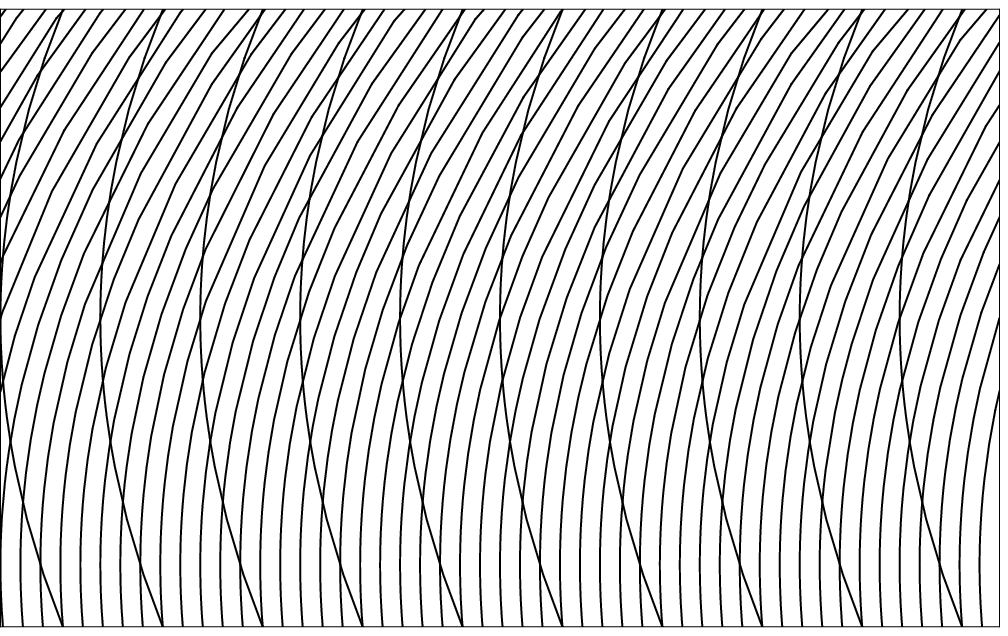}}
\caption{Period Class = (-0.6,-1.5)}
\label{fig:np3}
\end{figure}

\clearpage

%%%%%%%%%%%%%%%%%%%%%%%%%%%%%%%%%%%%%%%%%%%%%%%%%%%%%%%%%%%%%%%%%%%%%%%%
\section{Generalities on affine structures}

Let $\A^n$ denote real affine $n$-space with an orientation.
Let $\Aff(n)$ denote its group
of orientation-preserving affine automorphisms.
We use the standard linear representation.
$\A^n$ is the affine hyperplane $\R^n\times\{1\}\subset\R^{n+1}$
and $\Aff(n)$ the subgroup of $\GL(n+1;\R)$ consisting
of matrices
\begin{equation}\label{eq:standardlinear}
(A,b) := \bmatrix  A & b \\ 0 & 1 \endbmatrix 
\end{equation}
where $A\in\GL(n)$ and $b\in\R^n$. The affine action is defined by
\begin{equation*}
(v,1) 
%\xmapsto{(A,\ b )}
\, \stackrel{(A,b)}\longmapsto \,   
(A v + b, 1).
\end{equation*}

An {\em affine structure\/} on an oriented smooth manifold $M$ is a maximal
atlas of coordinate charts into $\A^n$ such that on
overlapping coordinate patches the transformation is locally 
the restriction of an element of $\Aff(n)$.
The manifold $M$ together with an affine structure is called
an {\em affine manifold.}  A local diffeomorphism
$N\stackrel{f}\longrightarrow M$ between affine manifolds is {\em
affine\/} if it is related to affine isomorphisms by the local affine
coordinate charts. Given a local diffeomorphism $f$, an affine
structure $M$ determines a unique affine structure on $N$ such that
$f$ is an affine map.  Thus covering spaces, and the universal
covering space $\tM\longrightarrow M$ in particular, of affine
manifolds are themselves affine manifolds. We shall assume for convenience
that our affine manifolds are path-connected and oriented.

Let $M$ be an affine manifold. 
The affine coordinate atlas on $M$ globalizes to a {\em developing map\/} 
$$ \dev : \tM \xrightarrow{} \A^n$$ 
from the universal covering $\tM$ into $\A^n$. This affine map
determines the affine structure on $M$, and is uniquely determined (for
a given affine structure) by its restriction to a coordinate patch. 
Conversely, any affine chart on a coordinate patch extends to a unique 
developing map. 

Let $M$ be an affine manifold and $p\in M$. A {\em affine germ\/} at $p$
is an equivalence class of affine diffeomorphisms 
\begin{equation*}
U\xrightarrow{\psi}\psi(U) \subset \A^n 
\end{equation*}
where $U$ is a neighborhood of $p$ in $M$ and $\psi(U)\subset\A^n$ is
open. 
Affine diffeomorphisms
$U\xrightarrow{\psi}\A^n$ and $U'\xrightarrow{\psi'} \A^n$
define the same germ if and only if they
agree on $U\cap U'$. Clearly $\Aff(n)$ acts simply transitively on the
set of affine germs at $p\in M$. Since an affine germ determines a
unique developing map, the developing map of a affine manifold is
unique up to composition with an affine transformation $g\in\Aff(n)$.

Deck transformations of $\tM\to M$ act by affine
automorphisms of $\tM$ and define the {\em affine holonomy
representation\/} $\pi_1(M)\stackrel{h}\longrightarrow\Aff(n)$ for
which $\dev$ is equivariant. The {\em developing pair\/} $(\dev,h)$ is
unique up to the action of $\Aff(n)$ where $\Aff(n)$ acts by 
left-composition on $\dev$ and by conjugation on $h$.

An affine manifold is {\em complete\/} if its developing map
is a diffeomorphism. Equivalently, an affine manifold is {\em
complete\/} if it is affinely diffeomorphic to a quotient space of
$\A^n$ by a discrete group $\Gamma$ of affine
transformations acting properly on $\A^n$.  
An {\em affine crystallographic group\/} is a discrete subgroup of
$\Aff(n)$ acting properly with compact quotient.
The classification of complete affine manifolds up
to affine diffeomorphism is equivalent to the classification of affine
crystallographic groups up to affine conjugacy.

Our goal is to construct a space whose points represent
inequivalent affine structures on a given closed manifold $S$.
Unfortunately such spaces are typically non-Hausdorff quotients
of singular spaces which are not smooth manifolds.
Thus we replace our goal with finding a smooth manifold
(or at least a smoothly stratified space)
with a group action, whose orbits represent affine diffeomorphism
classes of affine structures on $S$. 

To compare different affine manifolds, define a {\em marked 
affine manifold\/} as an orientation-preserving diffeomorphism
$S \xrightarrow{f} M$ where $M$ is an oriented affine manifold.  The
group $\Diff$ of all orientation-preserving diffeomorphisms
$S\xrightarrow{\phi}S$ acts by composition with $f$ 
on the set of  affine manifolds which are marked by $S$:
\begin{equation*}
S\xrightarrow{\phi}S \xrightarrow{f} M. 
\end{equation*}
Two affine manifolds $(f,M)$ and $(f',M')$ marked by $S$ are 
considered isomorphic if there exists an affine 
diffeomorphism $\varphi: M \rightarrow M'$ such that 
$\varphi \circ  f = f'$, they are   
\emph{isotopic} if there exists an affine 
diffeomorphism $\varphi: M \rightarrow M'$ such that 
$\varphi \circ  f$ and $f'$ are related by a diffeomorphism
in the identity component $\Diffo$ of $\Diff$. 

Let $x_0\in S$ be a basepoint, $\Pi:\tS\longrightarrow S$ the corresponding
universal covering space, and $\pi=\pi_1(S,x_0)$ its group of deck
transformations. Let $(f,M)$ be an affine manifold marked by $S$. 
An affine germ $[\psi]$ at 
$f(x_0)$ determines a developing pair $(\dev,h)$,
which is defined on $\tS$. This development pair depends only on the 
affine isomorphism class of the pair $\big((f,M),[\psi]\big)$. 

Let $\bigd$ denote the set of all affine isomorphism classes 
of pairs $\big((f,M),[\psi]\big)$. One can alternatively think 
of $\bigd$ as the set of developing pairs $(\dev,h)$. 
 Since the developing map $\dev$ 
determines the affine holonomy $h$, one can further regard
$\bigd$ as the set of developing maps of affine structures on $S$.
Affine holonomy defines a map
\begin{equation}\label{eq:hol}
\bigd \xrightarrow{\hol} \hpg
\end{equation}
which is invariant under the subgroup $\Diffob$ of
$\Diffo$ fixing $x_0$, since the elements of $\Diffob$ 
induce the identity map on $\pi = \pi_1(S,x_0)$.

These spaces are topologized as follows.
Give $\bigd$ the $C^{\infty}$-topology of developing maps: a
sequence of developing maps converges if and only if it and all of its
derivatives converge uniformly on compact subsets.
Give $\bigd/\Diffob$ the quotient topology.
If $\pi$ is finitely generated, then $\hpg$ is a real affine algebraic
set, and has the classical topology.
A sequence $h_m\in\hpg$ converges if and only if for all finite
subsets $\{\alpha_1,\dots,\alpha_N\}\subset\pi$, and all compact subsets
$K\subset\A^n$, the restrictions of $h_m(\alpha_i)$ to $K$ converge
uniformly for $i=1,\dots,N$.

In these topologies, $\hol$ is continuous.

Thurston~\cite{Thurston1} was the first to realize that $\hol$ defines a 
local homeomorphism
\begin{equation*}
\bigd/\Diffob \xrightarrow{\Hol}\hpg,
\end{equation*}
that is, $\hol$ is an open map and two nearby developing maps with
identical holonomy are isotopic by a basepoint-preserving isotopy.
See \cite{CaEpGr,Goldman2,Kapovich,Lok} for further discussion of
this general fact.

In general, the holonomy map $\Hol$ is not injective, even if $S$ is
compact. The simplest example occurs when $S=T^2$ and the holonomy has
cyclic image (covering spaces of Hopf manifolds).  See
Smillie\cite{Smillie}, Sullivan-Thurston~\cite{SuThu} and
Goldman~\cite{GoldmanF} for more exotic examples.

However, the restriction of $\Hol$ to {\em complete structures\/}
is injective and identifies complete structures with their
holonomy representations. 

Since $\pi$ is discrete, a representation $h\in\hpg$ defines a {\em 
proper action\/} if for all compact $K_1,K_2\subset\A^n$, the set
\begin{equation*}
\{\alpha \in \pi \mid \alpha(K_1)\cap K_2 \neq \emptyset \} 
\end{equation*}
is finite. If furthermore $h$ defines a {\em free action,\/}
then $\pi$ is torsionfree. In particular the quotient 
$\A^n/h(\pi)$ is Hausdorff, the quotient map
$\A^n\longrightarrow\A^n/h(\pi)$ is a covering space, and
$\A^n/h(\pi)$ inherits a smooth structure from $\A^n$.

Let $\hpgp$ denote the subset of $\hpg$ 
consisting of proper actions $h$ such that the
quotient $\A^n/h(\pi)$ is diffeomorphic to $S$. 
Let  $\DiffOb$ denote the subgroup of $\Diffb$ which 
induces the identity on $\pi_1(S, x_0)$. 
\begin{proposition}\label{prop:hol}
Let $\bigdc$ denote the subspace of\/ $\bigd$ corresponding
to complete affine structures.  Then the restriction
\begin{equation*}
\bigdc/\DiffOb \xrightarrow{\Hol}\hpgp
\end{equation*}
is a homeomorphism.
\end{proposition} 
We sketch the proof. Let $h\in\hpgp$.
A diffeomorphism $S \longrightarrow \A^n/h(\pi)$ defines
a complete structure on $S$ with holonomy $h$, showing surjectivity
of $\Hol$. If $S\xrightarrow{f} M$
is a marked affine structure on $M$ such that the 
holonomy $h$ is the
holonomy of a marked complete affine structure
\begin{equation*}
S\xrightarrow{f'} M' = \A^n/h(\pi),  
\end{equation*}
then the development map  
\begin{equation*}
\tS \xrightarrow{\tilde f} \tM  
\xrightarrow{\dev_M} \A^n
\end{equation*}
descends to a local diffeomorphism
$S \xrightarrow{F} M'$
inducing the isomorphism 
\begin{equation*}
\pi_1(S,x_0) \longrightarrow \pi_1(M',f(x_0)) = h(\pi)
\end{equation*}
defined by $h$.
Since $S$ is closed, the local
diffeomorphism is a covering space and since $\pi_1(F)$ is an
isomorphism, $F$ is a diffeomorphism. The composition 
$(f')^{-1}\circ F\in\Diff$ induces the identity on $\pi_1(S,x_0)$.
This completes the sketch.

Evidently $\hol$ is $\Aff(n)$-equivariant, where $\Aff(n)$ acts by
composition of the affine germ on $\bigd$ and by conjugation
on $\hpg$. The space $\bigdp$ of isomorphism classes of marked affine
manifolds $(f,M)$ identifies with the quotient 
$\bigd/\Aff(n)$, and $\hol$ defines a continuous map 
\begin{equation*}
\bigdp \longrightarrow \hpg/\Aff(n). 
\end{equation*}
Note that $\Aut(\pi_1(S,x_0))$ acts on $\hpg$ by right composition
on homomorphisms. Furthermore 
$\hol$ is equivariant with respect to the $\Diffb$-action on
$\bigd$ and the action on  $\hpg$ defined by the homomorphism
\begin{equation*}
\Diffb \longrightarrow \Aut(\pi_1(S,x_0)).
\end{equation*}
The space of isotopy classes of marked complete affine
manifolds $(f,M)$ identifies with the quotient 
$  \Def_c(S) = \bigdcp/\Diffob $
and is called the \emph{deformation space}.
The $\Diffb$-action on $\Def_c(S)$ 
factors through the discrete group
$\pi_0(\Diffb)$ and the  homeomorphism $\Hol$ of Proposition \ref{prop:hol} 
induces a surjective map
\begin{align*}
\Def_c(S) /&\pi_0(\Diffb) \longrightarrow  \\
&  \big(\hpgp/\Aff(n)\big)/\Aut(\pi_1(S,x_0)).
\end{align*}
The first space is the 
{\em moduli space,\/} consisting of affine diffeomorphism
classes of complete affine manifolds which are diffeomorphic 
to $S$. In this paper we show
that when $S= T^2$, then this space identifies naturally with
the quotient of $\R^2$ by $\SLtz$.

\subsection*{Example}
In general, neither $\bigdp$ nor 
$\hpg/\Aff(n)$ is Hausdorff.
Two-dimensional Hopf manifolds provide counterexamples.
Namely the quotients
\begin{equation*}
M_\epsilon := 
\bigg(\R^2 -\{(0,0)\} \bigg)\,/\,
\bigg\{ 
\bmatrix 2 & \epsilon \\ 0 & 2 \endbmatrix^n
\;\big|\;  n\in\Z \bigg\}
\end{equation*}
are affine 2-tori. All the structures with $\epsilon\neq 0$ are equivalent,
but are not equivalent to $M_0$. In particular the equivalence class
of $M_1$ contains the equivalence class of $M_0$ in its closure.

\subsection*{Remarks}
The conventions in this paper were chosen to simplify the discussion of
our specific example: complete affine structures on $T^2$. In general,
one may to want to work with markings which are only homeomorphisms or
homotopy equivalences, or with manifolds which are not oriented (or
nonorientable).  We have chosen to work in the oriented differentiable
category for convenience, thereby avoiding as many topological
technicalities as possible. Similarly one may want to avoid the discussion
of base-point preserving diffeomorphisms in discussing the mapping class group
(see \S 2), but in the interest of brevity we have chosen the present approach 
to avoid the general complications.

Also in the specific simple case $S=T^2$ which is treated 
in this paper, we have that $\hpgp$ coincides with the 
subset of $\hpg$ consisting
of \emph{all} proper actions of $\pi$, and the condition
on the quotient space is actually superfluous. 
In fact, if $\pi= \Z^2$ then $\A^2/h(\pi)$ must be 
diffeomorphic to the 2-torus from the classification 
of surfaces. (See \cite{BauesV} for further discussion,
and generalization to higher-dimensional examples.)
Moreover,  since $T^2$ is a surface,  
$\dif_0(T^2,x_0) = \dif^0(T^2,x_0)$, see \S 2.

\subsection*{Connections}
An affine structure is equivalent to a flat torsionfree affine connection.
Recall that an {\em affine connection\/} is a covariant differentiation
operation on vector fields
\begin{align*}
\Vect \times \Vect & \xrightarrow{\nabla} \Vect \\
(X,Y) & \longmapsto \nabla_X(Y)
\end{align*}
which is $\R$-bilinear, and for $f\in C^\infty(M)$,
satisfies
\begin{equation*}
\nabla_{fX} Y = f \nabla_{X} Y, \qquad
\nabla_{X}(f Y) = f \nabla_{X} Y + (Xf) Y
\end{equation*}
(where $Xf\in C^\infty(M)$ denotes the directional derivative of 
$f$ with respect to $X$).
The connection is torsionfree if and only if
\begin{equation*}
\nabla_X Y - \nabla_Y X = [X,Y] 
\end{equation*}
and flat if and only if
\begin{equation*}
\nabla_X \nabla_Y - \nabla_Y \nabla_X = \nabla_{[X,Y]}.
\end{equation*}
Conversely any such affine connection arises from an atlas of locally
affine coordinates as above. 
(Compare \cite{KobayashiNomizu}.)

Covariant differentiation extends to other tensor fields by
enforcing identities such as the following. For example, if
$\eta$ is a $1$-form, and $X\in \Vect$, then
$\nabla_X(\eta)$ is the $1$-form defined by:
\begin{equation*}
\big(\nabla_X(\eta)\big) (Y) = X \eta(Y) -\eta(\nabla_XY).
\end{equation*}
for all $Y\in\Vect$. If $\eta_1,\eta_2$ are $1$-forms, then
\begin{equation*}
\nabla_X(\eta_1\wedge \eta_2 ) =
\nabla_X(\eta_1)\wedge \eta_2 - \eta_1\wedge \nabla_X(\eta_2)
\end{equation*}
(and linearity) defines $\nabla_X\eta$ for any exterior $2$-form $\eta$.

A tensor field $\eta$ is {\em parallel\/} if $\nabla_X\eta = 0$
for all $Y\in\Vect$. Equivalently, in local affine coordinates
the coefficients of $\eta$ are constants. More generally $\eta$ 
is {\em polynomial of degree $\le r$\/} if and only if
\begin{equation*}
\nabla_{X_1} \cdots \nabla_{X_{r+1}} \big(\eta\big) = 0
\end{equation*}
for all $X_1,\dots,X_{r+1}\in\Vect$.

\section{The two-torus $T^2$}
Let $S$ be the two-torus $T^2$, that is the 
quotient space of $\R^2$ by the action of the
integer lattice $\Z^2 \subset \R^2$ by translations.
The quotient map 
\begin{equation*}
\tS = \R^2 \longrightarrow \R^2/\Z^2 = T^2
\end{equation*}
is a universal covering space, thereby providing global Euclidean coordinates
for $\tS$. Let the image of the origin $0\in\R^2$ be the basepoint
$x_0\in T^2$.

A developing map for $T^2$ is thus a local diffeomorphism 
\begin{equation*}
\R^2 \xrightarrow{\dev} \A^2 
\end{equation*}
which is equivariant with respect to the translation action of 
$\pi \cong \Z^2$ on $\R^2$, and an affine action defined by
a homomorphism $\rho: \pi \rightarrow \Aff(2)$:
\begin{equation*}
\pp + \bmatrix  m_1 \\ m_2 \endbmatrix \stackrel{\dev}\longmapsto 
(\gamma_1)^{m_1}(\gamma_2)^{m_2} \pp
\end{equation*}
where $\gamma_i = \rho(\alpha_i)$ for $i=1,2$ and
\begin{equation*}
\alpha_1 =  \bmatrix 1 \\ 0 \endbmatrix, \quad
\alpha_2 =  \bmatrix 0 \\ 1 \endbmatrix
\end{equation*}
generate the fundamental group $\pi\cong\Z^2$.

\subsection*{The modular group}
The {\em mapping class group\/} of $T^2$ is particularly tractable.
Every orientation-preserving
diffeomorphism $T^2\longrightarrow T^2$ is isotopic to
a basepoint-preserving diffeomorphism, and every isotopy between
basepoint-preserving diffeomorphisms can be deformed to a
basepoint-preserving isotopy. Furthermore the natural
map
\begin{equation*}
\Diffb \longrightarrow \Aut(\pi_1(S,x_0))
\end{equation*}
corresponds to an isomorphism
\begin{equation*}
\pi_0(\dif(T^2)) \xrightarrow{\cong} \Aut(\pi_1(T^2,x_0)) \cong \SLtz.
\end{equation*}
See \S6.4 of Stillwell~\cite{Stillwell} for discussion.

$\SLtz$ is explicitly realized by {\em linear\/} diffeomorphisms 
of $T^2$ as follows. A matrix $A\in\SLtz$ defines a linear automorphism
of $\R^2$ which normalizes the action of integer translations $\Z^2$.
Thus $A$ induces a diffeomorphism of the quotient $T^2 = \R^2/\Z^2$,
which preserves the basepoint $x_0$. Thus the extension
\begin{equation*}
1 \longrightarrow \dif_0(T^2) \longrightarrow \dif(T^2) \longrightarrow
\pi_0(\dif(T^2)) \cong \SLtz \longrightarrow 1
\end{equation*}
is split.

\section{Complete affine structures on $T^2$}\label{sec:structuresonT}

\subsection*{Euclidean structures}
Since translations of $\R^2$ are affine transformations,
$\Z^2$ acts as an affine crystallographic group and
$T^2$ inherits a complete affine structure from the structure
of $\R^2$ (its universal covering) as a vector space.
The generators $\gamma_1,\gamma_2$ are translations by linearly
independent vectors in $\R^2$ spanning a lattice in $\R^2$.
Furthermore any two lattices in $\R^2$ are equivalent by
a linear automorphism of $\R^2$, so different lattices 
in $\R^2$ define affinely diffeomorphic structures.
This affine structure is called the {\em Euclidean structure\/}
since it supports more refined geometric structures modelled
on Euclidean geometry. (Of course, most linear automorphisms
will not preserve this more refined metric structure.)

The affine equivalence class of the Euclidean structure 
will be the origin in the deformation
space. The mapping class group $\SLtz$
fixes this affine equivalence class.

\subsection*{Simply transitive affine actions}

However, $T^2$ carries many other non-equivalent affine structures.
These structures are not subordinate to Riemannian metrics, so we call
them {\em non-Riemannian.\/} All the structures arise from simply
transitive unipotent affine actions, and we briefly review the general
theory.

By Smillie~\cite{Smillie} or Fried-Goldman-Hirsch~\cite{FGH}, any
proper affine action of an abelian (or, more generally, nilpotent)
group is {\em unipotent.\/} That is, the linear parts are unipotent
linear transformations. Equivalently the corresponding elements of
$\GL(n+1)$ are unipotent.  Furthermore a subgroup consisting of
unipotent elements is conjugate to the upper triangular unipotent
subgroup of $\GL(n+1)$.

If $\Gamma\subset\Aff(n)$ is an affine crystallographic group which is
{\em nilpotent,\/} then there exists a unique simply connected {\em
unipotent\/} subgroup $G$ of $\Aff(n)$ containing $\Gamma$ 
which acts simply transitively
on $\A^n$. Indeed, $G$ is the {\em algebraic hull\/} of $\Gamma$ in
$\Aff(n)$, that is, the Zariski closure of $\Gamma$ in $\Aff(n)$. 

Let $\mathfrak{g}$ denote the Lie algebra of $G$. The exponential map
$\mathfrak{g}\xrightarrow{\exp} G$ is a diffeomorphism. 

For any $p\in\A^n$, the composition 
\begin{equation*}
\mathfrak{g}\xrightarrow{\exp} G \xrightarrow{\mathsf{ev}_p} \A^n
\end{equation*}
is a {\em polynomial diffeomorphism\/} of the vector space $\mathfrak{g}$
with the affine space $\A^n$, where.
\begin{align*}
G &\xrightarrow{\mathsf{ev}_p} \A^n \\
g &\longmapsto g(p)
\end{align*}
denotes evaluation at $p\in \A^n$.
Hence any two unipotent
affine crystallographic actions of a group $\Gamma$ are conjugate
by a polynomial diffeomorphism of $\A^n$ (\cite{FriedGoldman},
Theorem~1.21).

Furthermore this implies that the de Rham cohomology of the quotient $M =
\A^n/\Gamma$ is computable from the cohomology of polynomial exterior
differential forms on $M$ (\cite{Goldman3}).  Since $\Gamma$-invariant
polynomial tensor fields on $\A^n$ are invariant under the algebraic
hull $G$, every tensor at the origin $p\in\A^n$ extends uniquely to a
$\Gamma$-invariant polynomial tensor field on $\A^n$. Conversely, any
two $\Gamma$-invariant polynomial tensor fields which agree at $p$ are
equal. Thus the space of polynomial vector fields on $M$
identifies with the tangent space $T_p\A^n = \R^n$, the space of
polynomial $1-forms$ on $M$ identifies with the dual space $T_p^*\A^n
= (\R^n)^*$. Since a unipotent subgroup of $\Aff(n)$ preserves volume,
the top-dimensional cohomology is represented by a {\em parallel
volume form.\/}

Here is an explicit construction for $n=2$.
Every unipotent subgroup is conjugate to a subgroup of the multiplicatively
closed affine subspace $1 + {\mathfrak{N}_3}\subset \GL(3)$ where 
\begin{equation*}
{\mathfrak{N}_3} := \bigg\{ \bmatrix 0 & * & * \\ 0 & 0 & * \\  0 & 0 & 0 \endbmatrix\bigg\}
\end{equation*}
is the subalgebra of strictly upper triangular matrices. The algebra
${\mathfrak{N}_3}$ satisfies ${\mathfrak{N}_3}^3 = 0$, and thus the exponential and logarithm maps
are given by:
\begin{align*}
1 + {\mathfrak{N}_3} & \xrightarrow{\log} {\mathfrak{N}_3} \\
\gamma & \longmapsto (\gamma-1) - \frac12 (\gamma-1)^2 
\end{align*}
and
\begin{align*}
{\mathfrak{N}_3} & \xrightarrow{\exp} 1 + {\mathfrak{N}_3} \\
a & \longmapsto  1 + a + \frac12 a^2
\end{align*}
respectively. If $\gamma \in 1 +{\mathfrak{N}_3}$, then
\begin{align*}
\gamma^m &  = \exp \bigg( m \;\log(\gamma)\bigg)\\ &
 = 1\; +\; m\,\bigg( (\gamma-1) + \frac12  (\gamma-1)^2\bigg) 
\; +\; m^2\, \bigg(\frac12 (\gamma-1)^2\bigg).
\end{align*}
Taking $m$ to be real extends this to a definition of the unique
one-parameter subgroup in $1 + {\mathfrak{N}_3}$ containing the cyclic group
$\langle \gamma\rangle$.

In particular if $\gamma_1,\gamma_2\in 1 + {\mathfrak{N}_3}$ generate an action
of $\Z^2$, then
\begin{align}
\gamma_1^{m_1} \gamma_2^{m_2} & = 
1 + m_1\, \bigg( (\gamma_1-1) + \frac12 (\gamma_1-1)^2 \bigg) 
  + m_2\, \bigg( (\gamma_2-1) + \frac12 (\gamma_2-1)^2 \bigg)  \notag \\
& \qquad\qquad  + {m_1}^2\, \bigg( \frac12(\gamma_1-1)^2 \bigg)  
\; +\; m_1m_2\, \bigg( (\gamma_1-1)(\gamma_2-1) \bigg) \notag \\
& \qquad \qquad \qquad   + {m_2}^2\, \bigg( \frac12(\gamma_2-1)^2 \bigg) 
\label{eq:Malcev}
\end{align}
for $m_1,m_2\in\Z$. Taking $m_1,m_2\in\R$ defines the extension of a
unipotent $\Z^2$-action to a unipotent $\R^2$-action.  (This is a
special case of the {\em Mal'cev extension\/} of a finitely generated
torsionfree nilpotent group to a 1-connected nilpotent Lie group. See
Mal'cev~\cite{Malcev}, or Chapter II of Raghunathan~\cite{Raghunathan}
for the general construction.) For the application to complete affine
manifolds with virtually polycyclic fundamental group, see \S 1 of
\cite{FriedGoldman}.

\subsection*{Non-Riemannian structures}
Let $\epsilon\in\R$.
The mappings 
\begin{align*}
\A^2  & \xrightarrow{\Phi_\epsilon} \A^2  \\
\bmatrix x \\ y \endbmatrix
& \longmapsto \bmatrix x +\epsilon\, y^2/2 \\ y \endbmatrix
\end{align*}
define a one-parameter subgroup of {\em polynomial diffeomorphisms\/}
of $\A^2$. 

For $(s,t)\in\R^2$, let $T_{s,t}$ denote translation
by $(s,t)$:
\begin{align*}
\A^2  & \xrightarrow{T_{s,t}} \A^2  \\
\bmatrix x \\ y \endbmatrix
& \longmapsto \bmatrix x + s  \\ y + t \endbmatrix.
\end{align*}
$\Phi_\epsilon$ conjugates translations to the affine action
$\R^2 \xrightarrow{\rho_\epsilon} \Aff(2)$ defined by: 

\begin{equation}
\rho_{\epsilon}\bigg( \bmatrix s \\ t \endbmatrix \bigg) \,:\,
\bmatrix x \\ y \endbmatrix \,
\xmapsto{\Phi_\epsilon T_{s,t}\Phi_{-\epsilon}} 
% \big|\hspace{-5pt}\xrightarrow{~~\Phi_\epsilon T_{s,t}\Phi_{-\epsilon}~} 
\, \bmatrix x + \epsilon y t + (s + \epsilon t^2/2) \\ y + t  \endbmatrix.
\end{equation}

In the standard linear representation \eqref{eq:standardlinear},
this affine representation equals:
\begin{equation}\label{eq:defrhoe}
\rho_{\epsilon}\bigg( \bmatrix s \\ t \endbmatrix \bigg) \; = \; 
\bmatrix 1 & \epsilon\, t & s + \epsilon\, t^2/2 \\
0 & 1 & t \\ 0 & 0 & 1 \endbmatrix.
\end{equation}

%\begin{equation}
%\rho_{\epsilon}\bigg( \bmatrix s \\ t \endbmatrix \bigg) \; := \;
%\Phi_\epsilon T_{s,t}\Phi_{-\epsilon} \;:\;
%\bmatrix x \\ y \endbmatrix
%\longmapsto \bmatrix x + \epsilon y t + (s + \epsilon t^2/2) \\
%y + t  \endbmatrix.
%\end{equation}

\subsection*{Automorphisms}
The polynomial one-parameter group $\Phi_\epsilon$ on $\A^2$
interacts with three linear one-parameter groups on on $\A^2$:
\begin{itemize}
\item
{\em Dilations\/} defined by
\begin{equation}\label{eq:dilations}
\bmatrix x \\ y \endbmatrix \stackrel{\delta_\lambda}\longmapsto 
\bmatrix \lambda^2 x \\ \lambda y 
\endbmatrix 
\end{equation}
for $\lambda>0$; 
\item
{\em Shears\/} defined by
\begin{equation}\label{eq:shears}
\bmatrix x \\ y \endbmatrix \stackrel{\sigma_u}\longmapsto 
\bmatrix  x + u y  \\  y \endbmatrix 
\end{equation}
for $u\in\R$;
\item {\em Homotheties\/} defined by
\begin{equation}\label{eq:homotheties}
\bmatrix x \\ y \endbmatrix \stackrel{h_\lambda}\longmapsto 
\bmatrix  \lambda x   \\  \lambda y \endbmatrix 
\end{equation} for $\lambda \neq 0$.
\end{itemize}
Dilations and shears commute with $\Phi_\epsilon$, while homotheties
conjugate the various elements of the one-parameter group $\Phi$:
\begin{align}\label{eq:homot}
% (\delta_\lambda)^{-1} \Phi_{\epsilon} 
%\delta_\lambda & = \Phi_{\epsilon} \notag\\
%(\sigma_u)^{-1} \Phi_{\epsilon} \sigma_u & = \Phi_{\epsilon} \notag \\
(h_{\lambda})^{-1} \Phi_{\epsilon} h_{\lambda} & = \Phi_{\lambda\epsilon}\,  .
\end{align}
\eqref{eq:homot} implies that if $\epsilon\neq 0$, then the group
$\rho_\epsilon(\R^2)$ is affinely conjugate to $\rho_1(\R^2)$. 

When $\epsilon\neq 0$, different lattices in $\R^2$ give inequivalent
affine structures (unlike the Euclidean case $\epsilon = 0$).
The map $\Phi_\epsilon$ is a {\em polynomial developing map\/}
for this structure.

% The affine actions $\rho_\epsilon$ defined in \eqref{eq:defrhoe}

Every complete affine structure on $T^2$ arises from this construction.

Denote the group of orientation-preserving linear automorphisms of $\R^2$ 
by $\GLtrp$. Given
\begin{equation*}
R  = \bmatrix s_1 & s_2 \\ t_1 & t_2 \endbmatrix \in \GLtr 
\end{equation*}
and $\epsilon\in\R$, we define a marked affine manifold 
$(f,M)$  as follows. 
Represent $S=T^2$ as $\A^2/\pi$, where $\pi\cong\Z^2$ acts
by integer translations. Then $\Phi_\epsilon$ defines a developing
map 
\begin{equation} \label{eq:parameters}
\tS = \R^2 \xrightarrow{\tilde f= \Phi_\epsilon \circ R} \A^2
\end{equation}
and a holonomy homomorphism $h = h_{(\epsilon,R)}$
\begin{align*}
\pi &\xrightarrow{h} \Aff(2) \\
\bmatrix m_1 \\ m_2 \endbmatrix &\longmapsto
\rho_\epsilon\bigg(
\bmatrix m_1 s_1 + m_2 s_2 \\ m_1 t_1 + m_2 t_2\endbmatrix\bigg).
\end{align*}
The corresponding affine manifold $M$ is the quotient $\A^2/h(\pi)$
and the marking $f$ is the diffeomorphism induced by $\tilde f$ above.

The Mal'cev construction \eqref{eq:Malcev} extends every proper affine
action of $\Z^2$ on $\A^2$ to a simply transitive affine action
of $\R^2$ on $\A^2$.
By Kuiper \cite{Kuiper}, Furness-Arrowsmith \cite{FurArr}, 
Nagano-Yagi \cite{NaganoYagi}, Fried-Goldman~\cite{FriedGoldman},
this subgroup is conjugate to the subgroup
\begin{equation}\label{eq:go}
G_1 := \rho_1(\R^2) = 
\bigg\{
\bmatrix 1 & t & s + t^2/2 \\ 0 & 1 & t \\ 0 & 0 & 1 \endbmatrix \;
\bigg| \quad s,t\in\R \bigg\}.
\end{equation}
The affine representations $\rho_\epsilon$ define simply transitive
affine actions and the evaluation map
\begin{equation*}
G = \R^2  \longrightarrow \A^2
\end{equation*}
is the developing map for an {\em invariant affine structure\/} on the
abelian Lie group $G$. 

The automorphisms of the affine structure of $G$ 
which normalize the multiplications in $G$ are induced by the 
normalizer $\Nrm$ of the image of $G$ in $\Aff(2)$.
This normalizer is generated by two commuting one-parameter 
subgroups $\delta_\lambda$ of dilations and shears $\sigma_u$ 
defined in \eqref{eq:dilations}  and \eqref{eq:shears} respectively.

When $\epsilon = 0$, the corresponding normalizer is larger: in this case
the normalizer of the translation group is the full affine group.
Every automorphism of the vector group $\R^2$ is an element of $\GLtr$.

It follows that the two pairs 
$(\epsilon_i,R_i)$ (where $i=1,2$) determine the same point
in $\Def_c(T^2)$  if and only if:
\begin{itemize}
\item $\epsilon_1=\epsilon_2 = 0$; or
\item $\epsilon_i\neq 0$ for $i=1,2$ and
\begin{equation} \label{qs}
R_2 = \bigg(\frac{\epsilon_1}{\epsilon_2}\bigg)^2 \delta_\lambda \sigma_u R_1
\end{equation}
for some $\lambda>0$ and $u\in\R$.
\end{itemize}

\subsection*{Polynomial tensors}
An oriented non-Riemannian complete affine 2-torus $M$ 
admits a unique parallel area form $\omega_M$ of area $1$, as well
as a unique parallel $1$-form $\eta_M$ satisfying the conditions
of the following  Lemma.
We call $\eta_M$ the {\em canonical parallel $1$-form\/} on $M$.

\begin{lemma}\label{lem:conditions}
Let $M$ be a 2-torus with a non-Riemannian complete affine structure and
parallel 2-form $\omega_M$ such that
\begin{equation}\label{eq:area}
\int_M \omega_M = 1. 
\end{equation}
Then there exists a unique parallel vector field $\zeta_M$ such that
\begin{equation*}
\zeta_M = \nabla_{\xi_M}\xi_M 
\end{equation*}
for a polynomial vector field $\xi_M$ and
\begin{equation}\label{eq:wzx}
\omega_M(\zeta_M,\xi_M) = 1.
\end{equation}
There is a unique parallel $1$-form $\eta_M$ such that
\begin{equation}\label{eq:eta}
\eta_M = \iota_{\zeta_M}(\omega_M) 
\end{equation}
where $\iota$ denotes interior multiplication.
Furthermore $\xi_M$ is unique up to addition of a constant multiple of 
$\zeta_M$.
\end{lemma}
\begin{proof}
Choose coordinates so that $M= \A^2/\Gamma$ where $\Gamma\subset G_1$
is a lattice. Dilations distort area by
\begin{equation*}
(\delta_\lambda)^* dx\wedge dy = \lambda^3 dx\wedge dy,
\end{equation*}
so
replacing $\Gamma$ by $\delta_\lambda\Gamma \delta_\lambda^{-1}$, 
we may assume that
\begin{equation*}
\omega = dx\wedge dy 
\end{equation*}
induces an area form on $M$ with area $1$.

The $\Gamma$-invariant polynomial vector fields on $\A^2$ 
are $G_1$-invariant, and 
\begin{equation*}
\xi := \ddy + x \ddy, \qquad \zeta := \ddx
\end{equation*}
is a basis for the $\Gamma$-invariant polynomial vector fields on
$\A^2$. Every parallel vector field is a constant multiple of $\zeta$.
Furthermore $\xi,\zeta,\omega$ satisfy:
\begin{equation}\label{eq:conditions}
\nabla_\xi\xi = \zeta, \qquad \omega(\zeta,\xi) = 1. 
\end{equation}
Let $\xi_M,\zeta_M$ be the vector fields on $M$ induced by
$\xi,\zeta$ respectively. This establishes existence.

% Dilations and shears act on these forms by:
% \begin{align}\label{eq:dilact}
% (\delta_\lambda)^* \omega & = \lambda^3 \omega \notag \\
% (\delta_\lambda)^* \eta & = \lambda \eta \notag \\
% (\delta_\lambda)^* (\iota_\xi\omega) & = \lambda^2 (\iota_\xi\omega)
% \end{align}
% and
% \begin{align}\label{eq:shact}
% (\sigma_u)^* \omega & =  \omega \notag\\
% (\sigma_u)^* \eta & = \eta \notag \\
% (\sigma_u)^* (\iota_\xi\omega) & = \iota_\xi\omega + u \eta.
% \end{align}
Every polynomial vector field is of the form
\begin{equation}\label{eq:polyvec}
\xi' = a \xi +  b \zeta
\end{equation}
and determines a parallel vector field
\begin{equation*}
\zeta' = \nabla_{\xi'}\xi' = a^2\zeta
\end{equation*}
which satisfies
\begin{equation*}
\omega(\zeta',\xi') = a^3.  
\end{equation*}
Thus the coefficient in \eqref{eq:polyvec}
is $a=1$ and 
conditions \eqref{eq:conditions}
uniquely determine $\zeta' = \zeta$. 
Finally \eqref{eq:eta} uniquely determines $\eta_M$.
\end{proof}

\section{Periods as parameters}

Let $(f,M)$ be a marked complete affine 2-torus such that $M$ is 
non-Riemannian. Let $\eta_M$ be its canonical parallel $1$-form.
Let $\alpha_1,\alpha_2\in\pi$ be the standard basis for $\pi\cong\Z^2$. 
The pair
\begin{equation*}
\vp := \vp(f,M) := (t_1,t_2) \in \R^2
\end{equation*}
where
\begin{equation*}
t_i =  \int_{\alpha_i}  f^*\eta_M
\end{equation*}
for $i=1,2$,
is an invariant of $(f,M)$. 
Call $\vp(f,M)$ the {\em period class\/} of $(f,M)$.
If $M$ is Euclidean, define $\vp(f,M)$ to equal $(0,0)$.

\begin{lemma}\label{lem:nonexact}
If $M$ is non-Riemannian, the period class
$\vp(f,M)$ is nonzero.
\end{lemma}
\begin{proof}
Suppose that $\vp(f,M) = 0$.  Since the images of $\alpha_1,\alpha_2$
generate $H_1(M)$, the cohomology class $[\eta_M]\in H^1(M;\R)$ is
zero and $\eta_M$ is exact.  Thus $\eta_M=d\psi$ for a smooth
function $M\xrightarrow{\psi}\R$.  Since $\eta_M$ is nowhere
vanishing, $\psi$ is a submersion, contradicting compactness of $M$.
\end{proof}

Since the value of the period class depends only
on the cohomology class of $f^*\eta_M$, the map $\vp: [f,M] 
\mapsto \vp(f,M)$ is well defined on isotopy classes of marked tori. 

\begin{thm}\label{thm:periods} 
The period mapping
\begin{equation*}
\Def_c(T^2)  \xrightarrow{\vp} \R^2 
\end{equation*}
is a homeomorphism.
\end{thm}
\begin{proof}
%\bigskip\noindent{\bf computation of periods:}
Suppose that $(f,M)$ is a non-Riemannian marked torus, 
defined by parameters $(\epsilon,R) \in \R \times \GL_+(2,\R)$, $\epsilon \neq 0$,
as in $(\ref{eq:parameters})$. 
Since $\Phi_\epsilon$ is area-preserving $dx \wedge dy$ induces an area
form on $M$ which has total area $A = s_1t_2-s_2t_1$. Therefore $\omega_M$
is the area form on $M$ induced by
\begin{equation*}
\omega = A^{-1} dx\wedge dy.
\end{equation*}
A basis for the $G_\epsilon= \rho_\epsilon(\R^2)$-invariant polynomial 
vector fields on $\A^2$ is 
\begin{equation*}
\ddy + \epsilon y \ddx,\,  \ddx
\end{equation*}
and we compute $c\in\R$ for which 
\begin{equation*}
\xi := c\bigg(\ddy + \epsilon y \ddx\bigg)
\end{equation*}
induces polynomial vector fields $\xi_M,\zeta_M$
satisfying \eqref{eq:conditions}. Define
\begin{equation*}
\zeta  =  \nabla_\xi\xi  = c^2 \epsilon \ddx 
\end{equation*}
so that  $\omega(\xi,\zeta) = A^{-1} c^3 \epsilon$.
Thus \eqref{eq:wzx} implies 
\begin{equation*}
c =  \bigg(\frac{A}{\epsilon}\bigg)^{1/3}
\end{equation*}
so the canonical vector field, and  canonical parallel $1$-form 
on $M$ are represented by 
\begin{align*}
\zeta \;& =\; A^{2/3}\,\epsilon^{1/3}\, \ddx  \\
\eta \;=\; \iota_\zeta\omega \;& =\; 
\bigg(\frac{\epsilon}A\bigg)^{1/3} \, dy 
\end{align*}
respectively. The period class for $[f,M]$ is:
\begin{equation}\label{eq:periodformula}
\vp \;=\;  
\bigg(\frac{\epsilon}A\bigg)^{1/3}
\bmatrix  t_1 \\ t_2 \endbmatrix \; .
\end{equation}

\bigskip\noindent
\emph{$\vp$ is surjective:}
We construct an inverse $\vq$ to $\vp$ as follows.
Let $\vk = (k_1,k_2)\in\R^2$. Define
\begin{equation}\label{eq:defRt}
R_\vk := \Vert\vk\Vert^{2} 
\bmatrix k_2 & -k_1 \\ k_1 & k_2 \endbmatrix
\end{equation}
where 
\begin{equation*}
\Vert \vk \Vert = \sqrt{(k_1)^2 + (k_2)^2}
\end{equation*}
as usual. Then 
\begin{equation*}
\Psi_{\vk} := (R_\vk)^{-1} \circ\Phi_1 \circ R_\vk =
(Q_\vk)^{-1} \circ\Phi_{\Vert\vk\Vert^3} \circ Q_\vk
\end{equation*}
is a polynomial developing map,
where $Q_\vk$ is the orthogonal matrix
\begin{equation*}
Q_\vk := \Vert\vk\Vert^{-3}R_\vk = 
\bmatrix k_2/\Vert\vk\Vert 
& -k_1/\Vert\vk\Vert  
\\ k_1/\Vert\vk\Vert 
& k_2/\Vert\vk\Vert  \endbmatrix
\end{equation*}
when $\vk\neq 0$. (Compare \eqref{eq:homot}.) Explicitly,

\begin{equation}\label{eq:univdev}
\bmatrix x \\ y \endbmatrix \,
\xmapsto{\Psi_{\vk}}
\bmatrix x \\ y \endbmatrix  + \frac12 (k_1 x + k_2 y)^2 
\bmatrix k_2 \\ -k_1 \endbmatrix 
\end{equation}

so $\Psi_{\vk}$ depends \emph{smoothly} on $\vk\in\R^2$.
Thus in particular $\vk \mapsto \Psi_{\vk}$ defines a map 
\begin{equation*}
\R^2 \xrightarrow{\vq}  \Def_c(T^2)
\end{equation*}
which is  continuous. 

Note also that $\Psi_{\vk}^{-1} = \Psi_{-\vk}$.
Hence the corresponding holonomy representation 
\begin{equation}\label{eq:holhk} h_\vk = \Psi_{\vk} T  \Psi_{-\vk}
\end{equation}
depends smoothly on $\vk$. Furthermore $h_\vk$ equals 
the conjugate of $h_{(\epsilon,R)}$ by $R^{-1}$,  
where 
\begin{equation*}
\epsilon = \Vert\vk\Vert^3, \; R = Q_\vk.
\end{equation*} 
By \eqref{eq:periodformula}, the period class of $\vq(\vk)$ equals $\vk$.
Thus the composition
\begin{equation}\label{qi}
\R^2 \xrightarrow{\vq} \Def_c(T^2)  \xrightarrow{\vp} \R^2
\end{equation}
is the identity map on $\R^2$, and $\vp$ is surjective.

Note also that $\vq$ is surjective. Indeed,
since every element $R \in GL_+(2,\R)$ decomposes
as a product of an orthogonal and an upper-triangular matrix,
$(\ref{qs})$ implies that every element of $\Def_c(T^2)$ may be 
represented by a holonomy representation with parameters 
$(\epsilon, Q)$, where $Q \in SO(2)$ is
an orthogonal matrix.    
Therefore,
(\ref{qi}) shows that $\vq$ is a continuous bijection onto
$\Def_c(T^2)$.

\bigskip\noindent\emph{$\vp$ is continuous:~}
It follows from (\ref{eq:univdev}) that the set  
$\vS = \{\Psi_\vk \mid \vk \in \R^2\}$ is closed in 
the space of development maps $\DT$,
and hence $\Def_c(T^2)$ carries the quotient topology
from $\vS$. Since $\vq$ is a bijection, $\vS$ is actually
homeomorphic to $\Def_c(T^2)$. Write 
$\Psi_\vk = (\Psi^1_\vk, \Psi^2_\vk)$. Since 
\begin{equation*}
 k_1  = - \sqrt[3]{{\partial^2 \over {\partial y}^2} \Psi^2_\vk} \; \mbox{ and } \;
 k_2  =   \sqrt[3]{{\partial^2 \over {\partial x}^2} \Psi^1_\vk} \; ,
\end{equation*}  
$\vp: \vS \rightarrow \R^2$ is continuous.  
Therefore, $\vp$ is a continuous function on the 
deformation space as well. 
%
%Furthermore every point in
%$\Def_c(T^2)$ has an open neighborhood for which 
%$A = s_1t_2-s_2t_1$ is bounded away from zero. 
%
% Since every point has a neighborhood for which $A$ is bounded away from
% zero, $\vp$ is continuous.
%

\bigskip\noindent
\emph{$\vp$ is injective:}
Let $[(f,M)]$ be a point in
in $\Def_c(T^2)$. If $\vp(f,M) = 0$, then by Lemma~\ref{lem:nonexact},
$M$ must be Euclidean. All Euclidean structures correspond to a single
point in $\Def_c(T^2)$.

Now suppose that $M$ is non-Riemannian. Since $\vq$ is a bijection
onto $\Def_c(T^2)$, the composition (\ref{qi})
implies that $\vp$ is injective on non-Riemannian structures as well. 

It follows that $\vp$ is a homeomorphism with inverse $\vq$.
\end{proof}

\subsection*{The moduli space}
The moduli space of complete affine structures on 
$T^2$ is obtained as the quotient
space of $\Def_c(T^2)$ by the mapping class
group $\pi_0(\dif(T^2))$.

Explicitly, the action of the modular group $\SLtz$ on 
$\Def_c(T^2) = \R^2$ is given as follows.
Let $(f,M)$ be
an affine manifold marked by $S=T^2$, and $\phi_A$
the diffeomorphism of $T^2$ which corresponds to $A \in \SLtz$.
Then the action of $A$ on 
$\Def_c(T^2)$ is represented by $(f,M) \mapsto (f \circ \phi_A, M)$. 
\begin{lemma} $$ \vp(f \circ \phi_A, M) = A^t (\vp(f,M) )$$ 
\end{lemma}  
\begin{proof} We assume that $\tilde{f}: \A^2 \rightarrow \A^2$
is of the form $\Phi_\epsilon \circ R$, for some linear map
$R \in \GL(2,\R)$ with determinant $1$. By (\ref{eq:periodformula}), 
$\vp(f,M) = \epsilon^{1 \over 3} (t_1, t_2)^t$, where 
$(t_1, t_2)$ is the second row vector of $R$.  
Then $$\widetilde{f \circ \phi_A} = \tilde{f} \circ  A = \Phi_\epsilon \circ RA\; , $$
and $\vp(f \circ \phi_A, M)= \epsilon^{1 \over 3}  A^t (t_1, t_2)^t=
A^t (\vp(f,M))$.
\end{proof} 

Thus the moduli space identifies naturally with
the quotient of $\R^2$ by $\SLtz$. 

\subsection*{Differentiable families}
We  construct a \emph{differentiable family} of
complete affine tori over the deformation space: For $\vk \in \R^2$,
we let $h_\vk: \Z^2 \rightarrow \Aff(2)$ denote the holonomy
homomorphism for the developing map $\Psi_\vk$ 
(see proof of Theorem \ref{thm:periods}). We let $\Z^2$ act on
$\R^2 \times \A^2$ via
\begin{equation*}
\gamma: (\vk, v) \, \mapsto \, (\vk, h_\vk(\gamma) v ) \; . 
\end{equation*}
By \eqref{eq:holhk},
$\Z^2$ acts differentiably (in fact by polynomial
diffeomorphisms), and the quotient space $\mathcal{F}$ is a
differentiable manifold, and a torus-bundle over $\R^2$ via the
projection onto the first factor. This bundle is trivial as a
differentiable torus bundle.
The diffeomorphism  $(\vk, v) \mapsto (\vk, \Psi_\vk(v))$ 
induces an explicit trivialization  
$$  \R^2 \times T^2 \xrightarrow{}  \mathcal{F} \; . $$
The bundle 
$\mathcal{F} \rightarrow \R^2$ 
admits a system of locally trivializing charts in 
$\R^2 \times \A^2$ which define an affine structure on the fibers,
giving a differentiable family of affine manifolds over $\R^2$
(see \cite{KodairaSpencer} for the definition of differentiable family
of \emph{complex} manifolds). Each fibre of the corresponding 
differentiable family
\begin{equation*}
\mathcal{F} \longrightarrow  \Def_c(T^2)\;  
\end{equation*}
is a marked complete affine torus which represents its base-point in the
deformation space.

Not every differentiable family of complete affine tori
is induced by a \emph{smooth} map  
into the deformation space. Indeed, for $a \in \R$, let
us consider the development map $\Phi_a: \R^2 \rightarrow \A^2$ defined by
\[
\bmatrix x \\ y \endbmatrix \,
\xmapsto{\Phi_a} 
 \bmatrix x \\ y \endbmatrix + {1 \over 2}\, x^2 \bmatrix 0 \\ a \endbmatrix \; .
\]
Then $\Phi_a$, $a \in \R$,
constitutes a differentiable family of developing maps, and
there is a corresponding differentiable family of affine manifolds
$$ \mathcal{E} \longrightarrow \R $$
where each fibre  $\mathcal{E}_a$, $a \in \R$, is a complete 
affine torus marked by $\Phi_a$. (Since the parameter 
$a$ enters linearly, the family 
$ \mathcal{E} \rightarrow \R$ defines an affine fibration, 
see \cite{Tsemo}.) The corresponding curve
$ \R \rightarrow \Def_c(T^2)$
is a line in $\R^2$, parametrized by 
$$  a \mapsto \vp(\mathcal{E}_a) = (0, - a^{1/3}) \; . $$ 
It is not smooth in $a=0$. 

As mentioned in the introduction, such a phenomenon can not 
happen for deformations of elliptic curves and their induced
maps into $\h$, 
see \cite[Theorem 14.3]{KodairaSpencer}.

\end{document}